\documentclass[11pt,a4paper]{article}

\usepackage{t1enc}
\usepackage[latin1]{inputenc}

\usepackage{amsmath}
\usepackage{amssymb}
\usepackage{ulem}
\usepackage[left]{lineno}

\usepackage{colortbl}

\usepackage[T1]{fontenc} 
\usepackage{times}
\usepackage{graphicx}
\usepackage{setspace}
\usepackage{fancyheadings}
\newtheorem{Pro}{Proposition}[section]

\newtheorem{Cor}{Corollary}[section]
\newtheorem{Rem}{Remark}[section]

\newcommand{\ds}{\displaystyle}
\title{Stochastic perturbations and fisheries management}
\author{Patrice Loisel\thanks{Corresponding author, MISTEA, INRA,
    Montpellier SupAgro, 2 place Viala 34060 Montpellier, France,
    patrice.loisel@inra.fr, Tel: +33(0)4 99 61 29 04}}

\begin{document}

\maketitle
\begin{spacing}{1.62}

\begin{abstract}
  As most natural resources, fisheries are affected by random
  disturbances. The evolution of such resources may be modelled by a
  succession of deterministic process and random perturbations on
  biomass and/or growth rate at random times. We analyze the impact of
  the characteristics of the perturbations on the management of
  natural resources. We highlight the importance of using a dynamic
  programming approach in order to completely characterize the optimal
  solution, we also present the properties of the controlled model and
  give the behavior of the optimal harvest for specific jump kernels.
\end{abstract}

Recommendations for Resource Managers:
\begin{itemize}
\item In the context of updated biomass, for a centrally disturbed
  biomass and sufficiently high effort the optimal harvest increases
  with biomass jump rate
\item In the context of jointly updated biomass and growth rate, for a
  centrally disturbed biomass and sufficiently high effort the optimal
  harvest increases with the biomass jump rate
\item In the context of jointly updated biomass and growth rate, for
  sufficiently high effort the optimal harvest decreases with the
  growth jump rate
\end{itemize}

\section{Introduction}
The evolution of natural resources is most often disturbed by random
events. These disturbances occur at times that are not necessarily at
regular intervals. 
Hence, the management of natural resources must take into account the
characteristics of these disturbances.

The inclusion of stochastic perturbations in resource management has
been the subject of numerous articles in the literature, \cite
{Clark}, \cite {Pindyck}, \cite {Saphores}, \cite {Sethi}. Most of
these works concerns growth processes subject to perturbation
continuously or at regular intervals, whereas, as mentioned above,
perturbations occur most often at discrete random times. Our study
tries to take into account the random nature of not only the
perturbations magnitude of the resource but also the occurrence
times of these perturbations.

For systems with random perturbations 
the state variables are updated at random times and between these
random times, the state variables are governed by deterministic
processes. The most appropriate framework for the study of such
systems seems to be that of the Piecewise Deterministic Markov Process
(PDMP) \cite{Davis}. Using this framework, Hanson and Tuckwell
\cite{Hanson-Tuckwell} study the time to extinction of a population
with some specific growth function and random perturbation
structure. Applications in dynamic population are studied in
\cite{Ryan85}, \cite{Ryan}. Hanson \cite{Hanson} gives a panorama of
the models developed in various fields with this framework. We are
interested by the optimal management of fisheries in this
framework. The goal of this paper is to study the behavior of the
control variable with respect to the jump rate of the perturbation
process.

  In Section 2 we first present the resource growth model with its
  deterministic and stochastics components and the corresponding PDMP
  framework. Secondly, we express the controlled problem with updated
  biomass, we highlight the importance of using a dynamic programming
  approach in order to completely characterize the optimal solution,
  we also present the properties of the controlled model and give the
  application for a specific jump kernel. Finally in Section 4, we
  consider the case of updated biomass and growth rate.

\section{The model with updated biomass}

We assume that the evolution of the biomass in a fishery is mainly governed by a
determinist continuous process while it is observed or perturbated
only at discrete random times.
In absence of update, the evolution of the resource biomass $x(t)$ at
time $t$ is governed by a deterministic growth model:
\begin{align}
  {d x(t) \over dt} = \ & G(x(t))  - h(x(t),e(t)), \notag \\
\mbox{ with initial condition  }: x(0) = \ & x_0, \notag
\end{align}
with $0 < x_0 < K$. The parameter $K$ is the carrying capacity of the
studied system and $h$ is the harvest and $e(t)$ is the harvest rate.

We assume that $G$ is a differentiable concave growth function such
that $G(0)=0, G'(x) > 0$ for $ x < K$, $G'(x) < 0 $ for $x > K$. The
most common example is the logistic growth function
$\ds G(x)= rx (1-x/K)$ with the growth rate $r$.

We consider a resource submitted to random updates of the biomass
$\mathcal Y_1, \mathcal Y_2, ..$ at random times $\tau_1,
\tau_2,...$.
We assume that updates occur in a Poisson process i.e. that updates
occur independently of one another and randomly in time. The
distribution of times between successive updates is an exponential
distribution with mean $1/\lambda$:
\begin{align}
F(x) = & 1-e^{-\lambda x}, \notag 
\end{align}
with the constant jump rate $\lambda$. For each random time $\tau_i$,
the biomass is updated:
$$x(\tau_i^+) = \mathcal Y_i \overset{d}{\sim} \mathcal L (.|x(\tau_i)),
\mbox{ at time } \tau_i, \mbox{ for } i \geq 1.$$ where $\mathcal L$ is
a conditional distribution.

Hence the dynamics of the biomass can be described by a Piecewise
Deterministic Markov Process (PDMP) (\cite{Davis}). The random jump
process is described by the jump kernel operator $Q$. To each function
$\theta$
the operator $Q$ associates the function $Q[\theta]$ which is defined
by: $ Q[\theta] (x) = \int_{\mathcal Y} \theta(\mathcal Y) d\mathcal
L(\mathcal Y|x)$,
$Q[\theta] $ is assumed continuous and Lipschitz for all $\theta \in
\mathcal C(\Omega)$.

For example: $\mathcal Y_i = \mathcal Z_i x(\tau_i) $ with
$\mathcal Z_i \overset{d}{\sim} \mathcal U(\underline z,\overline z
)$, hence the associated jump kernel $Q$ is: \\
$\ds Q[\theta] (x) $
$ \ds = {1 \over \overline z-\underline z} \int_{\underline
  z}^{\overline z} \theta(zx) dz $
$ \ds = {1 \over (\overline z-\underline z)x} \int_{\underline
  zx}^{\overline zx} \theta(y) dy$.

\subsection{The biomass growth process}

In order to describe the biomass growth process, we define the
function $X(t;x,\tau)$ at time $t$. If the biomass was $x$ at time
$\tau$, the evolution of $X(t;x,\tau)$ is given by $(\mathcal S_{x_,\tau})$:
\begin{align}
  {d X(t;x,\tau) \over dt} = \ &
  g(X(t;x,\tau),e(t)) \equiv G(X(t;x,\tau))
-h(X(t;x,\tau),e(t)), \notag
\end{align} 
with initial condition: $X(\tau;x,\tau) = \ x$.\\
Knowing these characteristics, denoting $\tau_0$ at $x_0^+=x_0$, the
biomass growth process $\{ X_t(x_0): t \geq 0 \}$ starting with
biomass $x_0$ at initial time, may be expressed, for $i \geq 1$:
\begin{align}
                          & \ X_t = X(t;x_{i-1}^+,\tau_{i-1}), \ \tau_{i-1} < t \leq \tau_i, \notag \\
  \mbox{ while at time } \tau_i,& \ x_i^+ \overset{d}{\sim} \mathcal L (.|X_{\tau_i}).
                                  \notag
\end{align}
The process $\{ X_t(x_0): t \geq 0 \}$ starting at $x_0$ is composed
of the successive curves $X(t;x_0,0)$, $X(t;x_{\tau_1}^+,\tau_1)$,..,
$X(t;x_{\tau_i}^+,\tau_i)$, .... This process depends on the
successive random updated time $\tau_i$ and random jump at
corresponding $\tau_i$. \\
 Remark: by composite construction, for $ \tau_{i-1} < s < t < \tau_i$, we have:
 $X(t;X(s;x,\tau_{i-1}),s)=X(t;x,\tau_{i-1})$ and
 $X(t;x,\tau)=X(t-\tau;x,0)$.

\subsection{The control problem}

Given a biomass $x$ and an effort $e$, the instantaneous gain of
consumption $l(x,e)$ is determined. Therefore, we assume a regulator
maximizing expected discounted gain on an infinite horizon:
\begin{align}
\label{JE}
J(x_0,e(.)) = \ & E\Big[ \int_0^{+\infty} 
l(X_t(x_0),e(t)) e^{-\delta t} dt\Big],
\end{align}
with the $X_t(x_0)$ solution obtained with successive systems
$(\mathcal S_{x_0,0}), (\mathcal S_{x_{\tau_1}^+,\tau_1}), ...$.  The
expectation in Equation (\ref{JE}) is related to the successive random
updated time $\tau_i$ and random jump at corresponding $\tau_i$. The
effort $e$ is subject to the constraints:
$0 \leq e(t) \leq \overline e$ for all $t > 0$. The instantaneous gain
is assumed proportional to the effort: $l(x,e)=l_0(x) e$. Thus we
consider the function value $V$ defined by: $\ds V(x) = \max_{e(.) \in
  [0,\overline e]} J(x,e(.))$. \\
Assuming $V \in \mathcal C^1([0,K])$, we can formally deduce (see
Appendix A, with restrictive conditions \cite{Dempster}, \cite{Davis}
gives mathematical justification) the
Bellman Hamilton Jacobi (BHJ) equation:
\begin{align}
\label{BHJ}
\max_{e \in [0,\overline e]} [V'(x) g(x,e) - (\delta+\lambda) V(x) +
l(x,e) + \lambda Q[V](x)] & = 0.
\end{align}
The harvest is assumed proportional to the effort $e$ i.e.
$h(x,e)= h_0(x) e$, the BHJ equation becomes:
\begin{align}
\label{BHJm}
  \max_{e \in [0,\overline e]} [l_0(x)- h_0(x) V'(x)] e+ V'(x) G(x) -
  (\delta+\lambda) V(x) + \lambda Q[V](x) \ & = 0.
\end{align}
In Equation (\ref{BHJm}), the effort $e$ depends on the biomass
$X$. This highlights the existence of a critical value $x^*$ solution
of:
\begin{align}
\label{arc}
 l_0(x)-h_0(x) V'(x)=&\ 0.
\end{align}
As $0 \leq e \leq \overline e$, the optimal effort is a feedback
control $e^*(t) = \mathcal E(X_t)$ where the function $\mathcal E$ is defined by:
\begin{align}
\mathcal E(x) = \ &
\begin{cases}
0, & \mbox{ if }  l_0(x)- h_0(x) V'(x)< 0, \notag \\
{G(x) \over h_0(x)}, & \mbox{ if } l_0(x)- h_0(x) V'(x)= 0, \notag \\
\overline e, & \mbox{ if } l_0(x)- h_0(x) V'(x) > 0, \notag 
\end{cases}
\end{align}
where the function value $V$ is defined by: \begin{align}
\label{V} 
  [l_0(x)- h_0(x) V'(x)]_+ \overline e+ V'(x) G(x) - (\delta+\lambda)
  V(x) + \lambda Q[V](x) & = 0.
\end{align}
But Equation (\ref{arc}) is not sufficient to characterize the
critical values. By using a dynamic programming equation, we obtain a
complementary condition based on Euler-Lagrange condition.
 
\subsection{Euler-Lagrange condition}

The value function $V(x_0)$ is the solution to the optimization
problem:
$$V(x_0) = J(x_0,e^*(.))=\max_{e(.) \in [0,\overline e]}
E \Big[ \int_0^{+\infty} l(X_t(x_0),e(t)) e^{-\delta t} dt\Big],$$
with $X_t(x_0)$ solution obtained with successive systems
$(\mathcal S_{x_0,0}), (\mathcal S_{x_{\tau_1}^+,\tau_1}), ...$.

Using the strong Markov property, we may express \cite{Davis}:
\begin{align}
  V(x_0) = J(x_0,e^*(.))=& E_{\tau} \Big[\int_0^{\tau} l(X(t;x_0,0),e^*(t))
  e^{-\delta t} dt + Q[V](X(\tau;x_0,0)) e^{-\delta \tau}\Big]. \notag
\end{align}
We consider the first term in the right hand:
\begin{align}
  E_{\tau} \Big[\int_0^{\tau} l(X(t;x_0,0),e^*(t)) e^{-\delta t} dt\Big] = &
 \ \lambda \int_0^{+\infty} \int_0^{\tau} l(X(t;x_0,0),e^*(t)) e^{-\delta t} dt e^{-\lambda \tau} d\tau \notag \\
  = & \ \int_0^{+\infty} l(X(t;x_0,0),e^*(t)) e^{-(\delta+\lambda) t} dt,
  \notag
\end{align}
 when by inverting integration with respect to $t$ and $\tau$
Finally we obtain the dynamic programming equation:
$$ V(x_0) = \max_{e(.) \in [0,\overline e]} \int_0^{+\infty} [l(X(t;x_0,0),e(t))
+ \lambda Q[V](X(t;x_0,0))] e^{-(\delta+\lambda) t} dt.$$
To simplify the expressions, we denote $X(t,x_0) \equiv X(t;x_0,0)$.
We have the opportunity, as in the deterministic control case, to
deduce the expression of the effort $e(t)$ in terms of the biomass:
$$  l(X(t,x_0),e(t)) = l_0(X(t,x_0)) e(t) = {l_0 \over h_0}(X(t,x_0) (G(X(t,x_0))-\dot X(t,x_0))$$
and then we obtain the new form of the objective. Thus the
optimization problem becomes:
$$ V(x_0) =   \max_{X(.) \in \mathcal C_{x_0}} \int_0^{+\infty} \Big[{l_0 \over h_0}(X(t,x_0)) (G(X(t,x_0))-\dot X(t,x_0))
+ \lambda Q[V](X(t,x_0))\Big] e^{-(\delta+\lambda) t} dt. $$ $\mathcal C_{x_0}$ being
the set of admissible curves:
 $$ \mathcal C_{x_0} = \{ X(.) \in
 BC^1([0,K]), X(0)=x_0, G(X(t)) - h_0(X(t)) \overline e \leq \dot X(t)
 \leq G(X(t)) \},$$
 with $BC^1$ stands for the bounded with bounded derivative function
 defined on the interval $[0,K]$. We deduce:
\begin{Pro}
\label{Pro}
  Assuming that $V$ and $Q[V] \in B C^1([0,K])$, a critical value $x^*$
  is solution of the system of equations:
\begin{align}
 l_0(x)-h_0(x) V'(x)=&\ 0, \notag \\
\label{EL}
  (\delta+\lambda-G'(x))\Big[{l_0 \over h_0}\Big](x)= & \ \Big[{l_0 \over h_0}\Big]'(x) G(x) +\lambda [Q[V]]'(x), 
\end{align}
where the value fonction $V$ is solution of Equation (\ref{V}).
\end{Pro}
{\it Proof}: We have an implicit problem of Calculus of Variations.
 $X(.)$ stands for an interior solution, let $\mathcal
 L(.,.)$ the non actualized integrand: $\ds \mathcal I(X,\dot X) =
 {l_0 \over h_0}(X) (G(X)-\dot X) + \lambda Q[V](X)$, then $X(.)$ has to
 satisfy the Euler Lagrange condition:
$$\mathcal I_{X}(X(t),\dot X(t)) 
= {d \over dt} \mathcal I_{\dot X}(X(t),\dot X(t))- (\delta+\lambda)
\mathcal I_{\dot X}(X(t),\dot X(t)).$$ The Euler Lagrange condition
enhances:
\begin{align}
\Big[{l_0 \over h_0}\Big]'(X(t)) G(X(t) + {l_0 \over h_0}(X(t))
  G'(X(t))
  +\lambda [Q[V]]'(X(t))=  &  (\delta+\lambda) \Big[{l_0 \over h_0}\Big](X(t)).  \notag 
\end{align}
The differential equation is reduced to an algebraic Equation (\ref{EL}).
\hfill $\square$

Let $x^*(\lambda)$ the lower critical value, in order to avoid scaling
of function $V$ ($V$ is defined by Equation (\ref{V}) up to a
multiplicative constant for $x < x^*(\lambda)$), by using Equations
(\ref{arc}) and (\ref{EL}) becomes:
\begin{align}
\label{QVV}
(\delta-G'(x)+\lambda(1-{[Q[V]]'(x) \over V'(x)}))\Big[{l_0 \over h_0}\Big](x)=
& \Big[{l_0 \over h_0}\Big]'(x) G(x).
\end{align}
Similarly to the standard optimal control problem (without update) the
optimal effort $e$ is given by a function $\mathcal E$ of the biomass
$X$. But this function $\mathcal E$ depends on the jump rate $\lambda$ by the
intermediate of $x^*(\lambda)$.

\subsection{The value function and the optimal control}

We now analyze the behavior of the value function $V$ at critical value
$x^*$. Let $\ds A(x) = {l_0/h_0}(x)-V'(x)$. Then, for $x$ such
that $A(x) \neq 0$ we may define $\ds A'(x)=\Big[{l_0 \over h_0}\Big]'(x)-V''(x)$ and:
\begin{align}
\label{dV2}
 V''(x) & =\Big[{l_0 \over h_0}\Big]'(x)-A'(x),
\end{align}
so we deduce the regularity of the function value with respect to
biomass $x$:
\begin{Pro}
\label{Pro1}
The function value $V$ is continuously twice differentiable and:
\begin{align}
\label{arc2}
 V''(x^*) & =\Big[{l_0 \over h_0}\Big]'(x^*),
\end{align}
\end{Pro}
and with respect to jump rate $\lambda$:
\begin{Pro}
\label{Pro2}
For a sufficiently small value of jump rate $\lambda$, assuming
$l_0(x)=p q x -c$ and $h_0(x) = q x$ with price $p$, catchability $q$
and cost $c$, the critical value $x^*$ is an increasing (respectively
decreasing) function with respect to $\lambda$ if
$[Q[V]]'(x^*)- V'(x^*) > 0 $ (respectively $<0$). Moreover the
function value $V$ and $V_x'$ are continuously differentiable with
respect to jump rate $\lambda$.
\end{Pro} 
Proofs of Propositions are given in Appendix B. In the following
section, we will illustrate for a concrete case, with a specific jump
kernel, by a study of the sign of $[Q[V]]'(x^*) - V'(x^*)$. From
$V_x'$ continuously differentiable with respect to jump rate
$\lambda$, for a sufficiently small jump rate $\lambda$, Equations
(\ref{V}) and (\ref{EL}) has a unique solution $x^*$ so the function
$\mathcal E$ is given by:
\begin{align}
\mathcal E(x) = &
\begin{cases}
0, & \mbox{ if }  x < x^*, \notag \\
{G(x) \over h_0(x)}, & \mbox{ if } x=x^*, \notag \\
\overline e, & \mbox{ if } x > x^*, \notag 
\end{cases}
\end{align}
and finally:
\begin{Pro}
\label{Pro3}
  For a sufficiently small value of jump rate $\lambda$,
  assuming $l_0(x)=p q x -c$ and $h_0(x) = q x$ with price $p$,
  catchability $q$ and cost $c$, the value function is not three times
  differentiable, more precisely at $\lambda=0$:
\begin{align} 
  h_0(x^*) [x^2 V''(x)]'_+(x^*) = - {x^* \Sigma(x^*) \over \overline e-\mathcal E(x^*)}
  < & 0 < h_0(x^*) [x^2 V''(x)]'_-(x^*) = {x^* \Sigma(x^*) \over \mathcal E(x^*)}, \notag
\end{align} 
where $\ds \Sigma(x) = -x^2 V''(x)\Big[{G(x) \over x}\Big]' +(\delta-G'(x))
(V'(x)+V''(x) x)-G''(x) V'(x) x$ and $\Sigma(x^*) >0$. 
\end{Pro}
We now consider the growth function: $\ds G(x) = rx (1- {x/ K})$.
We may deduce the behavior of the critical value with respect to
growth rate $r$:
\begin{Pro}
\label{ProE}
For a sufficiently small value of jump rate $\lambda_x$ and $\lambda_r$,
assuming $l_0(x)=p q x -c$ and $h_0(x) = q x$, the critical value
$x^*(r)$ is an increasing function with respect to growth rate
$r$. 
\end{Pro} 
In Figure 1, we give an example of optimal evolution of the biomass
$x$. The dash line represents the level of the critical value $x^*$.

\includegraphics[scale=1,height=50mm,width=75mm,angle=0]{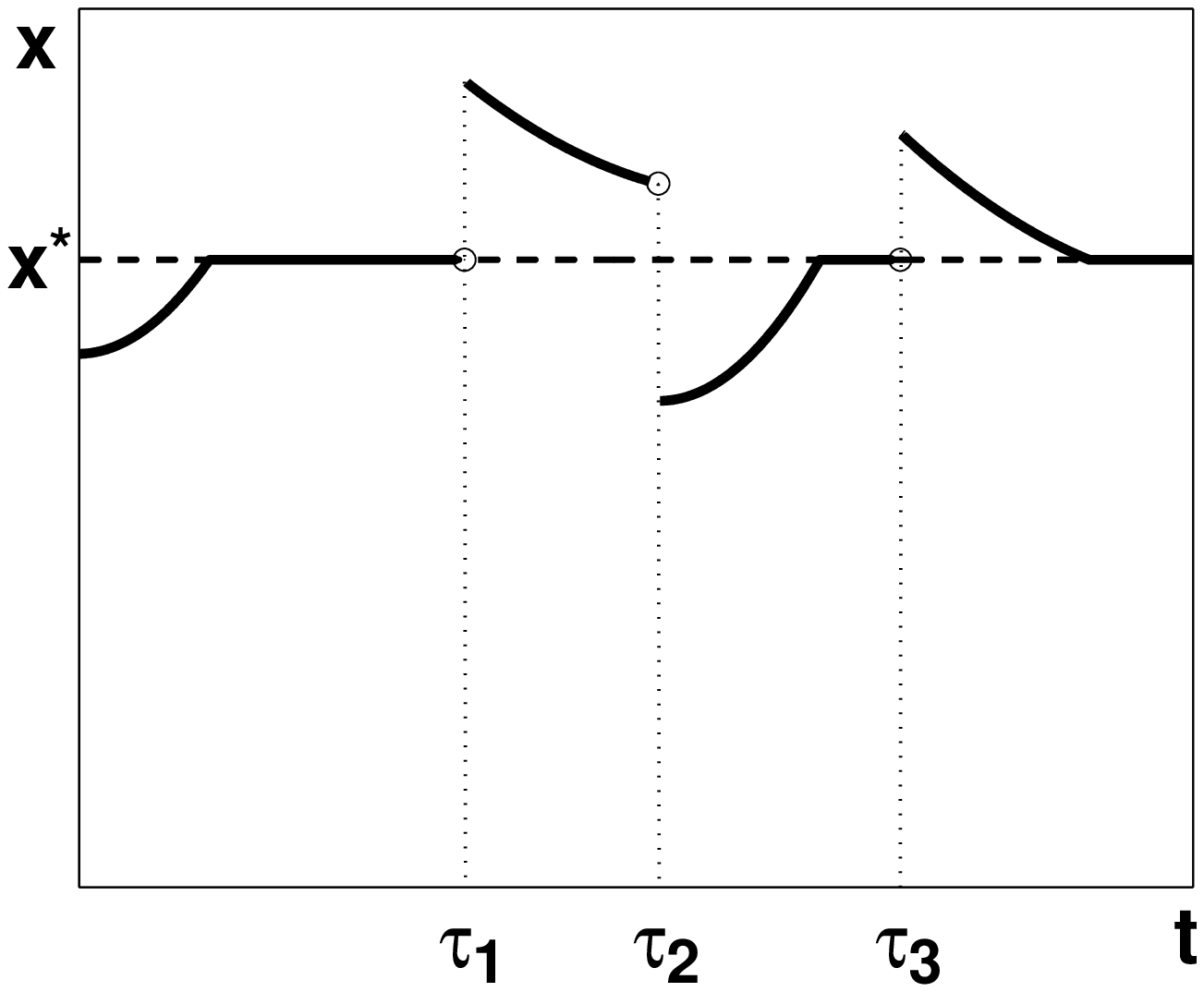}

Figure 1: Optimal evolution of biomass $x$ with biomass updated

\subsection{Application to specific jump kernel}

\begin{Pro}
\label{Pro4}
We assume
that 
the updated biomass $\mathcal Y_i$ is given by:
$\ds \mathcal Y_i = x(\tau_i)(1+ \mathcal Z_i \epsilon)$ where
$\ds \mathcal Z_i$ follows the distribution $\mathcal H$ (symmetric
and centered in $0$), for a sufficiently small value of jump rate
$\lambda_x$ and $x^*$ not close to $K$.

(i) if $E[\mathcal Z] \neq 0$ then the critical value is increasing
(respectively decreasing) with respect to jump rate $\lambda$ if $\ds
E[\mathcal Z]> 0$ (respectively $< 0$).

(ii) if $E[\mathcal Z] =0$ then the critical value is increasing
(respectively decreasing) with respect to jump rate $\lambda$ if $\mathcal E(x^*)$ is
smaller (respectively larger) than $\ds {\overline e \over 2}$.

More precisely in the latter case, with the growth function
$\ds G(x) = r x (1- {x / K})$, the critical value is increasing
(respectively decreasing) with respect to jump rate $\lambda$ if
$\ds x^* > K(1-{\overline e q / (2r)})$ (respectively $<$).

\end{Pro} 
Proofs are given in Appendix B. The given result for
$E[\mathcal Z] \neq 0$ is not surprising: for instance if
$E[\mathcal Z]>0$, higher jump rate leads to higher biomass, hence
higher possible harvest. If $E[\mathcal Z]=0$, the result is more
difficult to explain: for a sufficiently large value of $\overline e$,
it is optimal to use a higher level of critical value for the biomass.

\section{The model with updated biomass and growth rate}

We now consider a resource submitted two types of random updates
(biomass updates $\mathcal Y$, growth rate updates $\mathcal R$) at
random times: \\
- the times between two biomass updates
follows exponential distribution with mean $\ds {1 / \lambda_x}$\\
- the times between two growth rate updates follows
exponential distribution with mean $\ds {1 / \lambda_r}$. \\
Hence the random update time $\tau$ between two updates follows an
exponential distribution with mean
$\ds {1 /(\lambda_x+\lambda_r)}$. For each random time $\tau$, the
biomass or the growth rate is updated:
\begin{align}
x({\tau}^+) = \mathcal Y \overset{d}{\sim} \mathcal L_x (.|x({\tau}))
& \mbox{ with probability } {\lambda_x \over \lambda_x+\lambda_r} \notag \\
\mbox{ and }
  r({\tau}^+) = \mathcal R \overset{d}{\sim} \mathcal L_r (.|r({\tau}))
& \mbox{ with probability } {\lambda_r \over \lambda_x+\lambda_r}, \notag
\end{align}
where $\mathcal L_x$ and $\mathcal L_r$ are conditional distributions. \\
We assume that between two random updates the growth rate does not
change, i.e. the grow is piece-wise. The dynamics of the growth rate
is:
 $$\dot r(t) = 0.$$
 As in the previous case, the dynamics of the biomass can be described
 by a Piecewise Deterministic Markov Process (PDMP). The random jump
 process are described by the jump kernels $Q_x$ and $Q_r$. To each
 function $\theta$ of biomass $x$ and growth rate $r$, the functions
 $Q_x[\theta]$ and $Q_r[\theta]$ are defined by: $ Q_x[\theta] (x,r) =
 \int_{\mathcal Y} \theta(\mathcal Y,r) d\mathcal L_x(\mathcal Y|x),
 Q_r[\theta] (x,r) = \int_{\mathcal R} \theta(x,\mathcal R) d\mathcal
 L_r(\mathcal R|r)$.

\subsection{The biomass growth process}

In order to describe the biomass growth process we now define the
function $X(r,t;x,\tau)$ at time $t$. If the biomass was $x$ at time
$\tau$, the evolution of $X(r,t;x,\tau)$ is given by $(\mathcal S_{x_,r,\tau})$:
\begin{align}
  {d X(r,t;x,\tau) \over dt}
  = \ & G(X(r,t;x,\tau),r) -h(X(r,t;x,\tau),e(t)), \notag 
\end{align}
with initial condition: $X(r,\tau;x,\tau) = \ x$. \\
The process of the biomass $\{ X_t(x_0,r_0): t \geq 0 \}$ starting at
time $\tau_0=0$ with $x_0^+=x_0$, may be expressed for $i \geq 1$:
\begin{align}
 & \ X_t = X(r_{i-1},t;x_{i-1}^+,\tau_{i-1}), \ \tau_{i-1} < t \leq \tau_i, \notag \\
  \mbox{ where at time } \tau_i:& \  r_i=r_{i-1},\ x_i^+ \overset{d}{\sim} \mathcal L_x (.|X_{\tau_i}) \mbox{ with probability } {\lambda_x \over \lambda_x+\lambda_r} \notag \\ 
 \mbox{ and }           & \ x_i^+=X_{\tau_i},\  r_i \overset{d}{\sim} \mathcal L_r (.|r_{i-1}) \mbox{ with probability } {\lambda_r \over \lambda_x+\lambda_r}. \notag
\end{align}
Given a biomass $x$ and an effort $e$, we assume a regulator
maximizing expected discounted gain on an infinite horizon:
$$J(x_0,r_0,e(.)) = E\Big[ \int_0^{+\infty} 
l(X_t(x_0,r_0),e(t)) e^{-\delta t} dt\Big],$$ with $X_t(x_0,r_0)$ solution
obtained with successive systems $(\mathcal S_{x_0,r_0,0}),
(\mathcal S_{x_{\tau_1}^+,r_{\tau_1},\tau_1}), ...$.
Thus we consider the function value $V$ defined by: $\ds V(x,r) =
\max_{e(.) \in  [0,\overline e]} J(x,r,e(.))$. \\
Using the same formalism than for the updated biomass model and
assuming $V \in \mathcal C^1([0,K]\star[\underline r,\overline r])$,
denoting $\lambda.Q[V] = \lambda_x Q_x[V]+ \lambda_r Q_r[V]$, we can
formally deduce the corresponding Bellman Hamilton Jacobi (BHJ)
Equation:
\begin{align}
\label{BHJr}
\max_{e \in [0,\overline e]} [V_x'(x,r) g(x,r,e) -
(\delta+\lambda_x+\lambda_r) V(x,r) + l(x,e) + \lambda. Q[V](x,r)] & = 0.
\end{align}
The BHJ Equation becomes:
\begin{align}
  \max_{e \in [0,\overline e]} [l_0(x)- h_0(x) V_x'(x,r)] e+ V_x'(x,r) G(x,r) -
  (\delta+\lambda_x+\lambda_r) V(x,r) + \lambda. Q[V](x,r) \ & = 0. \notag
\end{align}
In this Equation, the effort $e$ depends on the biomass $x$ and the
growth rate $r$. As $0 \leq e \leq \overline e$, the optimal effort is
a feedback control $e^*(t) = \mathcal E(X_t,R_t)$ where the function $\mathcal E$ is defined by:
\begin{align}
\mathcal E(x,r) = &
\begin{cases}
0, & \mbox{ if }  l_0(x)- h_0(x) V_x'(x,r)< 0, \notag \\
{G(x,r) \over h_0(x)}, & \mbox{ if } l_0(x)- h_0(x) V_x'(x,r)= 0, \notag \\
\overline e, & \mbox{ if } l_0(x)- h_0(x) V_x'(x,r) > 0. \notag 
\end{cases}
\end{align}
Hence: \begin{align}
\label{Vr} 
  [l_0(x)- h_0(x) V_x'(x,r)]_+ \overline e+ V_x'(x,r) G(x,r) - (\delta+\lambda_x
+\lambda_r)  V(x,r) + \lambda. Q[V](x,r) & = 0.
\end{align}
The critical value $x^*(r)$ is the solution of the equation:
\begin{align}
\label{arcr}
 l_0(x)-h_0(x) V_x'(x,r)=&\ 0.
\end{align}
But this equation is not sufficient to characterize the critical
value. By using a dynamic programming equation, we obtain a
complementary condition based on Euler-Lagrange condition.

The value function $V(x_0,r_0)$ is the solution to the optimization
problem:
$$V(x_0,r_0) = J(x_0,r_0,e^*(.))=\max_{e(.) \in [0,\overline e]}
E \Big[ \int_0^{+\infty} l(X_t(x_0,r_0),e(t)) e^{-\delta t} dt\Big],$$ with
$X_t(x_0,r_0)$ solution of the system $(\mathcal S_{x_0,r_0,0})$. \\
Using the same reasoning than in the previous section, we obtain the dynamic programming equation:
$$ V(x_0,r_0) = \max_{e(.) \in [0,\overline e]} \int_0^{+\infty} [l(X(r_0,t;x_0,0),e(t))
+ \lambda. Q[V](X(r_0,t;x_0,0),r_0)] e^{-(\delta+\lambda_x+\lambda_r) t} dt.$$
\begin{Pro}
\label{Pror}
  Assuming that $V, Q_x[V]$ and $Q_r[V] \in B C^1([0,K]\star[\underline r,
  \overline r])$, a critical value $x^*(r)$ is solution of the system of
  equations:
\begin{align}
 l_0(x)-h_0(x) V_x'(x,r)=&\ 0 \notag \\
\label{ELr}
  \mbox{ and }  (\delta+\lambda_x+\lambda_r-G_x'(x,r))\Big[{l_0 \over h_0}\Big](x)= &  \Big[{l_0 \over h_0}\Big]'(x) G(x,r) +[\lambda .Q[V]]'(x), 
\end{align}
where the value function $V$ is solution of the Equation (\ref{Vr}).
\end{Pro}
For fixed $r$, using the same reasoning than in the previous section:
\begin{align}
  (\delta-G_x'(x,r)+\lambda_x(1-{[Q_x[V]]'_x(x,r) \over V'_x(x,r)})
  +\lambda_r(1-{[Q_r[V]]'_x(x,r) \over V'_x(x,r)}))\Big[{l_0 \over h_0}\Big](x)= &  \Big[{l_0 \over h_0}\Big]'(x) G(x,r) \notag
\end{align}
and replacing $\lambda$ by $\lambda_x$ (respectively $\lambda_r$) and
$\lambda Q[V]$ by $\lambda_x Q_x[V]$ (respectively $\lambda_r Q_r[V]$) we can
obtain the equivalent of the Propositions \ref{Pro1} and \ref{Pro2}.
\begin{Pro}
\label{Proxr}
The function value $V$ is twice differentiable with respect to biomass
$x$ and:
\begin{align}
\label{Vxx2}
V''_{x^2}(x^*(r),r) = & \Big[{l_0 \over h_0}\Big]'(x^*(r)).
\end{align}
Assuming $l_0(x)=p q x -c$ and $h_0(x) = q x$ with price $p$,
catchability $q$ and cost $c$, the critical value $x^*(r)$ is an
increasing (respectively decreasing) function with respect to biomass
jump rate $\lambda_x$ if $[Q_x[V]]_x'(x^*(r),r) -V_x'(x^*(r),r) > 0$
(respectively $<0$) and is an increasing (respectively decreasing)
function with respect to growth jump rate $\lambda_r$ if
$[Q_r[V]]_x'(x^*(r),r)-V_x'(x^*(r),r) > 0$ (respectively
$<0$). Moreover the function value $V$ and $V_x'$ are continuously
differentiable with respect to jump rate $\lambda_x$, and jump growth
rate $\lambda_r$.
\end{Pro} 
In the following section, we will illustrate for a concrete case, with
a specific jump kernels for biomass (respectively growth rate), by a
study of the sign of $[Q_x[V]]'_x(x^*,r) - V'_x(x^*,r)$ (respectively
$[Q_r[V]]'_x(x^*,r) - V'_x(x^*,r)$). For each growth rate $r$ and for
sufficiently small values of jump rate $\lambda_x, \lambda_r$,
Equations (\ref{arcr}) and (\ref{ELr}) have a unique solution
$x^*(r)$, the function $\mathcal E$ is given by:
\begin{align}
\mathcal E(x,r) = &
\begin{cases}
0, & \mbox{ if }  x < x^*(r), \notag \\
{G(x,r) \over h_0(x)}, & \mbox{ if } x=x^*(r), \notag \\
\overline e, & \mbox{ if } x > x^*(r), \notag 
\end{cases}
\end{align}
\begin{Rem}
  The expression of $\mathcal E(x,r)$ has the same formulation than in the case
  without updates but due to the difference between the Propositions
  \ref{Pro} and \ref{Pror}, the corresponding $x^*(r)$ differs.
\end{Rem}
\begin{Pro}
\label{LemB}
The value function $V$ is twice differentiable but not third
differentiable and at $\lambda=0$:
\begin{align}
 \label{Vxr}
 V''_{xr}(x^*(r),r) & = 0, \\
\label{Vxxr}
{V'''_{x^2r}}^-(x^*(r),r) = - {\Sigma_1(x^*(r),r) \over \mathcal E(x^*(r),r)} & <
0 < {V'''_{x^2r}}^+(x^*(r),r) = {\Sigma_1(x^*(r),r) \over \overline e
  -\mathcal E(x^*(r),r)}
\end{align}
and where $\ds \Sigma_1(x,r) = {\delta \over r}{l_0 \over
  h_0^2}(x)$.
The third derivatives of the function value $V$ with respect to
biomass $x$ and growth rate $r$ are linked by:
\begin{align}
\label{dVxr}
{V'''_{x^2r}}^\pm(x^*(r),r) {x^*}'(r) + {V'''_{xr^2}}^\pm(x^*(r),r)& =0. 
\end{align}
\end{Pro}
Proofs are given in Appendix C.  In Figure 2, we give an example of
optimal evolution of the biomass $x$ with alternatively biomass and
growth rate updating. The dash lines represent the successive levels
of the critical value $x^*(r)$,

\includegraphics[scale=1,height=50mm,width=75mm,angle=0]{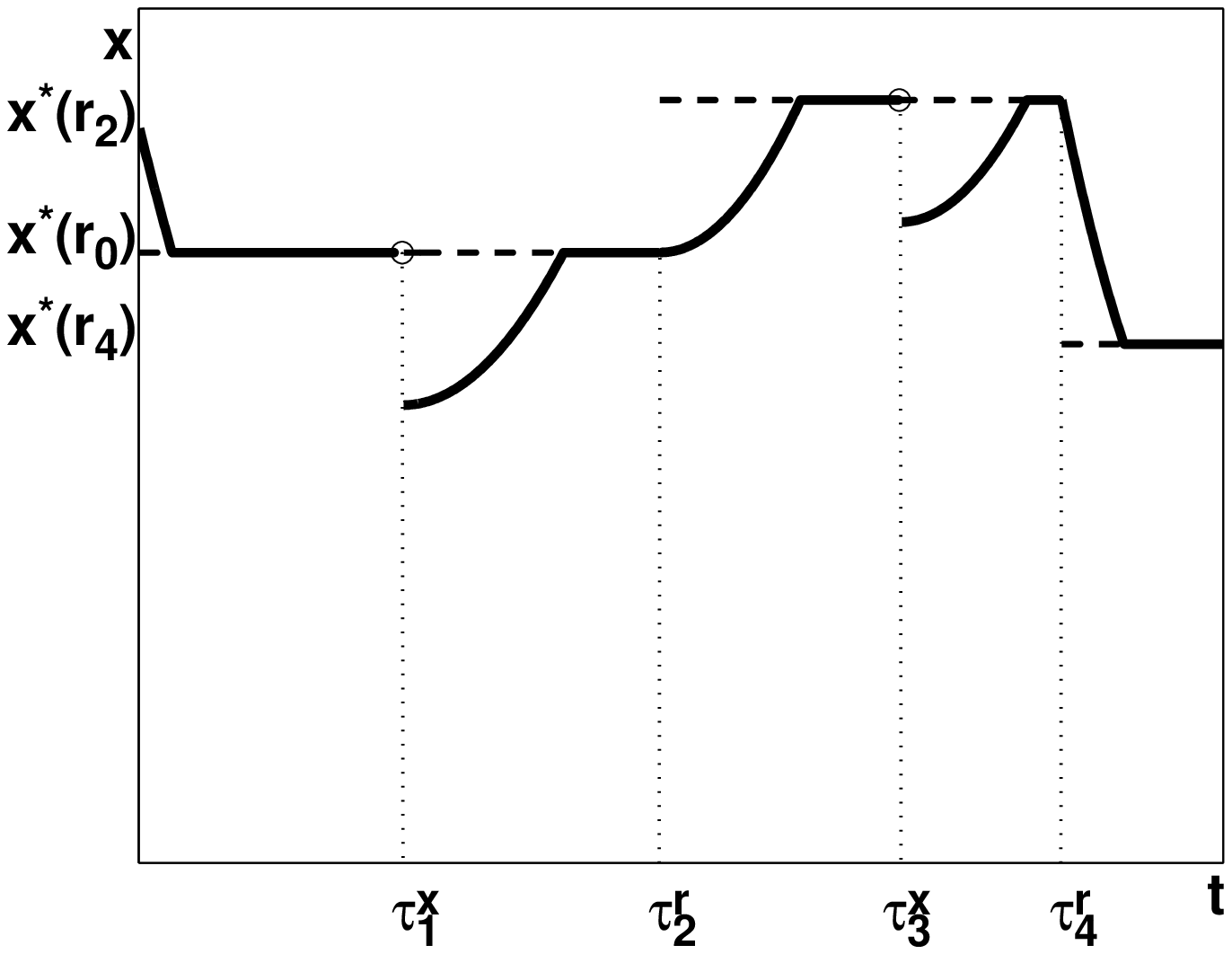}

Figure 2: Optimal evolution of biomass $x$ with biomass and growth rate updated

\subsection{Application to specific jump kernels}

\begin{Pro}
\label{ProF}
Assuming the updated biomass $\mathcal Y_i$ (respectively, the updated
growth rate $\mathcal R_i$) given by:
$\ds \mathcal Y_i = x({\tau_i^x})(1+ \mathcal Z_i \epsilon)$
(respectively, $\mathcal R_i = r({\tau_i^r}) (1+ \mathcal S_i \xi)$) where
$\ds \mathcal Z_i$ (respectively, $\ds \mathcal S_i$) follows the distribution
$\mathcal H_x$
(respectively, $\mathcal H_r$). \\
For a sufficiently small value of biomass jump rate $\lambda_x$ and
$x^*(r)$ not close to $K$:

(i) if $E[\mathcal Z] \neq 0$ then the critical value $x^*(r)$ is
increasing (respectively, decreasing) with respect to biomass jump rate $\lambda_x$ if
$\ds E[\mathcal Z]> 0$ (respectively, $< 0$).\\
(ii) if $E[\mathcal Z] =0$ then the critical value $x^*(r)$ is
increasing (respectively, decreasing) with respect to biomass jump rate $\lambda_x$ if
$\mathcal E(x^*(r),r))$ is larger (respectively, smaller) than $\ds {\overline e \over
  2}$.

For a sufficiently small value of growth jump rate $\lambda_r$ and $r$
not
close to $\underline r$ and $\overline r$: \\
(iii) the critical value $x^*(r)$ is increasing (respectively, decreasing)
with respect to growth jump rate $\lambda_r$ if $\mathcal E(x^*(r),r)$
is larger (respectively, smaller) than $\ds {\overline e \over 2}$.
\end{Pro}
We deduce the following properties:
\begin{Cor}
  (i) The behavior of the critical value with respect to biomass jump
  rate $\lambda_x$ is of the same type that in the previous case with
  update of the biomass. \\
  (ii) If $E[\mathcal Z] = 0$, the behavior of the critical value
  $x^*(r)$ is reversed with respect the biomass jump rate $\lambda_x$
  and the growth jump rate $\lambda_r$. \\
  (iii) The behavior of the critical value $x^*(r)$ with respect to
  growth jump rate $\lambda_r$ is independent of the expectation
  $E[\xi]$.
\end{Cor}

\section{Conclusion}
In this article, we consider the evolution of a fishery following a
continuous process and submitted to random updates at random times, we
present the appropriate PDMP framework. We express the control problem
with biomass updates, we highlight the importance of using a dynamic
programming approach in order to completely characterize the critical
value of the control. We give conditions which permit to deduce the
behavior of the optimal control effort with respect to jump rate. An
application to a specific jump kernel shows the possible variety of
behavior of the optimal effort with respect to the random
structure. For a centrally disturbed biomass and sufficiently high
effort the optimal harvest increases with biomass jump rate Finally we
study the more complex case for which biomass and growth rate in the
dynamics are updated. For a centrally disturbed biomass and
sufficiently high effort the optimal harvest increases with the
biomass jump rate and decreases with the growth jump rate.

\section*{Appendix A}

Let the value function defined by: $\ds V(x) = \max_{e(.) \in
  [0,\overline e]} J(x,e(.))$.  Consider a time $t > 0$, by the strong
Markov property, the criteria satisfies:
\begin{align}
V(x)=  J(x,e^*(.)) = \ & E_{\tau} \Big[(\int_0^{\tau}
  l(X(u,x),e^*(u))e^{-\delta u} du +
  Q[V](X(\tau,x)) e^{-\delta \tau}) I_{\tau < t}  \notag \\
  & + (\int_0^t l(X(u,x),e^*(u))e^{-\delta u} du + V(X(t,x))e^{-\delta
    t}) I_{t < \tau}\Big], \notag
\end{align}
\begin{align}
  E_{\tau} \Big[\int_0^{\tau} l(X(u,x),e^*(u)) e^{-\delta u} du I_{\tau <
    t}\Big] = \ &
  \lambda \int_0^{t} \int_0^{\tau} l(X(u,x),e^*(u)) e^{-\delta u} du \ e^{-\lambda \tau} d\tau \notag \\
  & \mbox{ then by inverting integration with respect to $t$ and $\tau$ } \notag \\
  = \ & \int_0^{t} l(X(u,x),e^*(u)) (e^{-(\delta+\lambda) u}-e^{-\lambda
    t -\delta u}) du, \notag
\end{align}
and hence, rearranging the terms:
\begin{align}
  V(x) = &\max_{e(.) \in [0,\overline e]} \Big[\int_0^{t}
  [l(X(\tau,x),e(\tau)) + \lambda Q[V](X(\tau,x))] e^{-(\delta+\lambda)
    \tau} d\tau + e^{-(\delta+\lambda) t} V(X(t,x))\Big], \notag
\end{align}
where $V(x)$ is independent of $t$, hence, formally differentiating with
respect to $t$:
$$\max_{e \in [0,\overline e]} [V'(x) g(x,e) - (\delta+\lambda) V(x) +
l(x,e) + \lambda Q[V](x)] = 0.$$

\section*{Appendix B}

\noindent {\bf Proof of Proposition \ref{Pro1}}. In $\Omega -
\{x|A(x)=0\}$, the value function is smooth and is solution of:
\begin{align}
  \eta A(x) h_0(x)\overline e + V'(x) G(x)
  -(\delta+\lambda)V(x)+\lambda Q[V](x) = & \ 0, \notag
\end{align}
with $\eta=0$ if $A(x) < 0$ and $\eta = 1$ if $A(x) > 0$.

By differentiation:
\begin{align}
\label{d2V}
  \eta (A'(x) h_0(x)+A(x) h'_0(x))\overline e + V''(x) G(x)
  -(\delta+\lambda-G'(x))V'(x)+\lambda [Q[V]]'(x) = &\ 0. 
\end{align}

Let $A_-'(x^*)$ and $A_+'(x^*)$ (respectively, $V_-''(x^*)$ and $V_+''(x^*)$)
the left and right derivatives of $A$ (resp $V'$) at the critical
value $x^*$ (if $A(x) (x-x^*) > 0$ in the vicinity of $x^*$ and the
reverse if not). From Equations (\ref{EL}) and (\ref{arc}), using
$A(x^*)=0$, at the critical value $x^*$:
\begin{align}
  V_-''(x^*) G(x^*) - \Big[{l_0 \over h_0}\Big]'(x^*) G(x^*)= & 0 \notag \\
 \mbox{ and } A_+'(x^*) h_0(x^*)\overline e + V_+''(x^*) G(x^*) - \Big[{l_0 \over h_0}\Big]'(x^*) G(x^*)= & 0, \notag 
\end{align}
hence respectively using the Equation (\ref{dV2}):
\begin{align} 
- A_-'(x^*) G(x^*) & = \ 0 \notag \\
\mbox{ and } A_+'(x^*) (h_0(x^*)\overline e -G(x^*)) & = \ 0, \notag
\end{align}
then $A_-'(x^*)= A_+'(x^*)=0$, so $A$ is differentiable 
and $V$ is twice countinuously differentiable. \hfill $\square$

\noindent {\bf Proof of Proposition \ref{Pro2}}. In order to determine
the behavior of the critical value $x^*$ in the vicinity of
$\lambda=0$, we differentiate the Equation (\ref{QVV}) with respect to
jump rate $\lambda$ to obtain:
\begin{align} 
  \Big[-G''(x^*) {l_0 \over h_0}(x^*) + (\delta-2 G'(x^*))\Big[{l_0 \over
    h_0}\Big]'(x^*)-\Big[{l_0 \over h_0}\Big]''(x^*)G(x^*)\Big]{x^*}_\lambda'(0)= & {l_0 \over
    h_0}(x^*) ({[Q[V]]'_x \over V'_x}(0,x^*)-1). \notag
\end{align}
Using expression of $l_0(x)$ and $h_0(x)$ the second equation becomes:
\begin{align} \Big[-G''(x^*) {l_0 \over h_0}(x^*) +(\delta +2({G(x^*) \over
    x^*}-G'(x^*)))\Big[{l_0 \over h_0}\Big]'(x^*)\Big]{x^*}_\lambda'(0)= & {l_0 \over
    h_0}(x^*) ({[Q[V]]'_x \over V'_x}(0,x^*)-1). \notag
\end{align}
From $G''(x) < 0$ and $G(0)=0$, we deduce that $G(x)-xG'(x)\geq 0$ and
the left term in the brackets is positive hence ${x^*_\lambda}'(0)$ and
$[Q[V]]'(x)-V'(x)$ for $\lambda=0$ has the same sign.\\
We derive Equations (\ref{V}) for $x < x^*(\lambda)$ and (\ref{arc})
with respect to jump rate $\lambda$ at $\lambda=0$, by using Equation
(\ref{d2V}):
\begin{align}
 V''_{x \lambda}(0,x^*) G(x^*) + Q[V](0,x^*)-V(0,x^*)= & \ \delta V'_\lambda(0,x^*) \notag \\
\mbox{ and } \Big[{l_0 \over h_0}\Big]'(x^*) {x^*}'_\lambda(0) = & \ V''_{x \lambda}(0,x^*). \notag
\end{align}
Hence: 
$V$ and $V_x'$ are continuously differentiable with respect to jump
rate $\lambda$.
\hfill $\square$

\noindent {\bf Proof of Proposition \ref{Pro3}}. By differentiation of
Equation (\ref{d2V}) multiplied by $x$, for $\lambda=0$:
\begin{align}
  \eta [x(A'(x) h_0(x)+A(x)h'_0(x))]' \overline e 
+ [x^2 V''(x)]' {G(x) \over x}  - \Sigma(x)&  =  \ 0. \notag
\end{align}
Hence due to $A(x^*)=A'(x^*)=0$:
\begin{align}
  [x^2 V''(x)]'_-(x^*) G(x^*) = & x^* \Sigma(x^*) \notag \\
\mbox{ and }  [x^2 V''(x)]'_+(x^*) (G(x^*)-h_0(x^*) \overline e) + (x^* \Big[x \Big[{l_0
  \over h_0}\Big]'(x)h_0(x)\Big]' (x^*)- {x^*}^3 V''(x^*) \Big[{h_0 \over
  x}\Big]'(x^*))\overline e= & x^* \Sigma(x^*). \notag
\end{align}
From expression of $l_0$ and $h_0$, $\ds \Big[x \Big[{l_0 \over
  h_0}\Big]'(x)h_0(x)\Big]'(x)=0$ and $\ds \Big[{h_0 \over x}\Big]'(x)=0$ so:
\begin{align}
  [x^2 V''(x)]'_+(x^*) (G(x^*)-h_0(x^*) \overline e) = & x^* \Sigma(x^*). \notag
\end{align}
From expression of $\mathcal E(x^*)$:
\begin{align}
  h_0(x^*) [x^2 V''(x)]'_-(x^*) \mathcal E(x^*) = & x^* \Sigma(x^*), \notag \\
  h_0(x^*) [x^2 V''(x)]'_+(x^*) (\mathcal E(x^*)-\overline e) = & x^* \Sigma(x^*). \notag
\end{align}
From concavity of $G$, $G(x)-xG'(x) \geq 0$ and $G''(x) < 0$, \\
hence $\Sigma(x) >0$, so $\ds [x^2 V''(x)]'_-(x^*) >0$ and
$\ds [x^2 V''(x)]'_+(x^*) < 0$. \hfill $\square$

\noindent {\bf Proof of Proposition \ref{ProE}}. In order to determine
the behaviour of an optimal critical value $x^*$, we differentiate the
Equation (\ref{EL}) with respect to growth rate $r$ at
$\lambda=0$ to obtain:
\begin{align}
  \Big[-G''_{x^2}(x^*,r) {l_0 \over h_0}(x^*) + (\delta-2 G'_x(x^*,r))\Big[{l_0
    \over h_0}\Big]'(x^*)-\Big[{l_0 \over h_0}\Big]''(x^*)G(x^*,r)\Big]{x^*}'(r)= &
  {\delta \over r} {l_0 \over h_0}(x), \notag \\
  \Big[-G''_{x^2}(x^*,r) {l_0 \over h_0}(x^*) +(\delta +2({G(x^*,r) \over
    x^*}-G'_x(x^*,r)))\Big[{l_0 \over h_0}\Big]'(x^*)\Big]{x^*}'(r)= & {\delta
    \over r} {l_0 \over h_0}(x). \notag
\end{align}
From $G''_{x^2}(x,r) < 0$ and $G(0)=0$, we deduce $G(x,r)-x G'_x(x,r)
\geq 0$ and the left term in the brackets is positive, hence from
positive marginal gain. \hfill $\square$

\noindent {\bf Proof of Proposition \ref{Pro4}}. For biomass $x$ not
close to $K$: $\ds Q[V](x) = \int_X V(x(1+z \epsilon)) dH(z)$, \\
so: $\ds [Q[V]]'(x)-V'(x)= \int_X [(1+z \epsilon) V'(x(1+z
\epsilon))-V'(x)] dH(z)$,
\begin{align}
  \mbox{where } V'(x(1+z\epsilon))(1+z\epsilon) -V'(x) = &\ z\epsilon (x V''(x)+V'(x)) \notag \\
                                                         & +\int_0^{z \epsilon} (z\epsilon-t) [V''((1+t)x)(1+t)x+V'((1+t)x)]'_t dt \notag \\
  = &\ z \epsilon (x V''(x)+V'(x)) \notag \\
                                                         & + \int_0^{z \epsilon} (z\epsilon-t) [V'''((1+t)x)(1+t)x^2+2V''((1+t)x)x] dt. \notag
\end{align}
For a sufficiently small $\epsilon$, for all $x \neq x^*$ the integral
can be approximated by $\ds{z^2 \epsilon^2 \over 2} (x^2 V'''(x)+2x
V''(x))$ (i.e. $\ds {z^2 \epsilon^2 \over 2} [x^2 V''(x)]'$) in the
vicinity of critical biomass $x^*$. Hence at critical biomass $x^*$:

(i) $\ds \lim_{\epsilon \rightarrow 0} {[Q[V]]'(x^*)-V'(x^*) \over
  \epsilon}= E[\mathcal Z](x^* V''(x^*)+V'(x^*))= E[\mathcal Z]\Big[x {l_0
  \over h_0}\Big]'(x^*)$,

(ii) $\ds \lim_{\epsilon \rightarrow 0} {[Q[V]]'(x^*)-V'(x^*) \over
  \epsilon^2}= {E[\mathcal Z^2] \over 4} ([x^2 V''(x)]_-'(x^*) +[x^2
V''(x)]'_+(x^*))$.
 
From Proposition \ref{Pro3}, 
$\ds [x^2 V''(x)]_-'(x^*) +[x^2 V''(x)]'_+(x^*)= {x^* \Sigma(x^*)
  \over h_0(x^*)} {\overline e- 2 \mathcal E(x^*) \over \overline e- \mathcal E(x^*)}$,
hence using Proposition \ref{Pro2}, $x'(0) > 0$ if
$\ds \mathcal E(x^*) < {\overline e \over 2}$ and $x'(0) < 0$ if
$\ds \mathcal E(x^*) > {\overline e \over 2}$. \hfill $\square$

\section*{Appendix C}

\noindent {\bf Proof of Proposition \ref{LemB}}. From total
differentiation of Equation (\ref{arcr}) with respect to growth rate
$r$:
$$(l_0'(x^*(r))-h_0'(x^*(r)) V''_{xx}(x^*(r),r)){x^*}'(r) = h_0(x^*(r)) V''_{xr}(x^*(r),r),$$
and using Equation (\ref{Vxx2}) we deduce Equation (\ref{Vxr}).
Let $\ds A(x,r) = {l_0 \over h_0}(x) -V_x'(x,r)$. \\
In $\Omega -\{x|A(x,r)=0\}$ the value function is smooth and is
solution of:
\begin{align}
  \label{d2Vr}
\eta A(x,r) h_0(x) \overline e + V_x'(x,r) G(x,r) - (\delta+\lambda_x
+\lambda_r)  V(x,r) + \lambda. Q[V](x,r) & = 0,
 \end{align}
 with $\eta = 0$ if $A(x,r) < 0$ (i.e. $x < x^*(r)$) and $\eta=1$
 if $A(x,r) > 0$ (i.e. $x > x^*(r)$). \\
 Using Equation (\ref{Vxx2}), from differentiation of Equation
 (\ref{d2Vr}) with respect to biomass $x$ and growth rate $r$ at
 $\lambda_x=\lambda_r=0$, at optimal critical value $x^*(r)$:
\begin{align}
  (G(x^*(r),r)-\eta h_0(x^*(r)) \overline e)V'''_{x^2r}(x^*(r),r) &+
  V''_{x^2}(x^*(r),r)
  G'_r(x^*(r),r) + V'_x(x^*(r),r) G''_{rx}(x^*(r),r) =0, \notag \\
  (G(x^*(r),r)-\eta h_0(x^*(r)) \overline e)V'''_{x^2r}(x^*(r),r) &+
  {1 \over r} (\Big[{l_0
    \over h_0}\Big]'(x^*(r)) G(x^*(r),r) + {l_0 \over h_0}(x^*(r)) G'_x(x^*(r),r)) =0, \notag \\
  \mbox{ then using Equation (\ref{ELr}) } & \notag \\
  (G(x^*(r),r)-\eta h_0(x^*(r)) \overline e)V'''_{x^2r}(x^*(r),r) &+
  {\delta \over r} {l_0    \over h_0}(x^*(r)) =0, \notag \\
  h_0(x^*(r))(\mathcal E(x^*(r),r)- \eta \overline e)V'''_{x^2r}(x^*(r),r) &+
  {\delta \over r} {l_0 \over h_0}(x^*(r)) =0.
\end{align}
Hence, from positive marginal gain we deduce the expression of
$\ds {V'''_{x^2r}}^-(x^*(r))$ and $ {V'''_{x^2r}}^+(x^*(r))$. \\
From total differentiation of Equation (\ref{Vxr}) with respect to
growth rate $r$ we deduce (\ref{dVxr}). \hfill $\square$

\noindent {\bf Proof of Proposition \ref{ProF}} (i) and (ii) Using the
same reasoning than in the previous section and replacing $\lambda$ by
$\lambda_x$ and $\lambda Q[V]$ by $\lambda_x Q_x[V]$ we obtain the
same result than in the Proposition \ref{Pro4}.\\
(iii) For growth rate $r$ not close to $\underline r$ and
$\overline r$:
\begin{align}
  Q_r[V](x,r)-V(x,r) & = \int_{\underline r}^{\overline r} (V(x,r(1+s
\xi))-V(x,r)) dH_r(s), \notag \\
\mbox{hence: } [Q_r[V]]'_x(x,r)-V'_x(x,r)&= \int_{\underline r}^{\overline r}[V'_x(x,r(1+s \xi))-V'_x(x,r)] dH_r(s) \notag 
\end{align}
and:
\begin{align}
  V'_x(x,r(1+s \xi))-V'_x(x,r)& =s\xi r V''_{xr}(x,r)+\int_0^{s
    \xi} (s \xi-t) r^2 V'''_{xr^2}(x,r(1+t)) dt. \notag
\end{align}
For a sufficiently small $\xi$, the integral can be approximated by
$\ds {s^2 \xi^2 \over 2} r^2 ({V'''_{xr^2}})_{\{-sign(s)\}}(x,r)$ in
the vicinity of the critical biomass $x^*(r)$ for all growth rate
$r$. Hence, due to Equation (\ref{Vxr}) :
$\ds [Q_r[V]]'_x(x,r)-V'_x(x,r)= {E[\mathcal R^2] \xi^2 \over 4} r^2
({V'''_{xr^2}}^-(x^*,r)+{V'''_{xr^2}}^+(x^*,r))+ O(\xi^3) $.
From 
Equations (\ref{Vxxr}) and (\ref{dVxr}),
$\ds {V'''_{xr^2}}^-(x^*,r)+{V'''_{xr^2}}^+(x^*,r)$ is proportional
with the same sign to $\overline e -2 \mathcal E(x^*(r),r)$. \hfill
$\square$

\section{Figure Legends}

Figure 1: Optimal evolution of biomass $x$ with biomass updated

\noindent Figure 2: Optimal evolution of biomass $x$ with biomass and growth rate updated

\end{spacing}
\end{document}